\documentclass[11pt]{amsart}%
\usepackage{amsfonts}
\usepackage{amsmath}
\usepackage{amssymb}
\usepackage{graphicx}%
\providecommand{\U}[1]{\protect\rule{.1in}{.1in}}
\newtheorem{theorem}{Theorem}
\theoremstyle{plain}

\numberwithin{equation}{section}
\begin{document}
\title[On the hypercontractivity of the Bohnenblust--Hille inequality]{A note on the hypercontractivity of the polynomial Bohnenblust--Hille inequality}
\author{Daniel Pellegrino }
\address{Departamento de Matem\'{a}tica, UFPB, Jo\~{a}o Pessoa, PB, Brazil}
\email{dmpellegrino@gmail.com and pellegrino@pq.cnpq.br}
\thanks{2010 Mathematics Subject Classification. 32A70, 32A50, 46G25, 47H60}
\keywords{Sidon constants, homogeneous polynomials, Bohnenblust--Hille inequality}

\begin{abstract}
For $\mathbb{K}=\mathbb{R}$ or $\mathbb{C}$ and $m$ a positive integer, we
remark that there is a constant $C$ so that, for all $r\in\lbrack1,\frac
{2m}{m+1}],$ the supremum of the ratio between the $\ell_{r}$ norm of the
coefficients of any nonzero $m$-homogeneous polynomial $P:\ell_{\infty}%
^{n}\left(  \mathbb{K}\right)  \rightarrow\mathbb{K}$ and its supremum norm is
dominated by $C^{m}\cdot n^{\left(  \frac{m}{r}-\frac{m+1}{2}\right)  }$ and,
moreover, we prove that the exponent $\frac{m}{r}-\frac{m+1}{2}$ is optimal. 

\end{abstract}
\maketitle

%\dedicatory{Submitted to the JFA: Category -- Short Communication}

\section{Introduction}

If $P:\ell_{\infty}^{n}\left(  \mathbb{C}\right)  \rightarrow\mathbb{C}$ is an
$m$-homogeneous polynomial defined by%
\[
P(z)={\sum\limits_{\left\vert \alpha\right\vert =m}}a_{\alpha}z^{\alpha}%
\]
and $\mathbb{D}^{n}$ is the unit polydisc in $\mathbb{C}^{n},$ let%
\[
\left\Vert P\right\Vert _{r}:=\left(  {\sum\limits_{\left\vert \alpha
\right\vert =m}}\left\vert a_{\alpha}\right\vert ^{r}\right)  ^{1/r}\text{ and
}\left\Vert P\right\Vert _{\infty}:=\sup_{z\in\mathbb{D}^{n}}\left\vert
P(z)\right\vert .
\]
The Sidon constant $S_{1}(m,n)$ (see \cite{annals, seip}) is the smallest
constant satisfying the inequality
\begin{equation}
\left\Vert P\right\Vert _{1}\leq S_{1}(m,n)\left\Vert P\right\Vert _{\infty}
\label{sid}%
\end{equation}
for all $m$-homogeneous polynomials $P:\ell_{\infty}^{n}\left(  \mathbb{C}%
\right)  \rightarrow\mathbb{C}.$ From \cite{annals, seip} (see also the
references therein) we know that there is an absolute constant $C_{\mathbb{C}%
}>0$ such that
\begin{equation}
S_{1}(m,n)\leq C_{\mathbb{C}}^{m}\cdot n^{\frac{m-1}{2}}. \label{comp}%
\end{equation}
It is also known that this result is sharp (we refer to \cite{annals, seip}
and the references therein).

The inequality (\ref{sid}) is related to the famous polynomial
Bohnenblust--Hille inequality for complex scalars (see \cite{bh}). Since the
polynomial Bohnenblust--Hille inequality is hypercontractive for both real and
complex scalars (see \cite{CMPS, annals}), there exists a constant
$C_{\mathbb{K},BH}>1$ such that%
\begin{equation}
\left\Vert P\right\Vert _{\frac{2m}{m+1}}\leq C_{\mathbb{K},BH}^{m}\left\Vert
P\right\Vert _{\infty} \label{09}%
\end{equation}
for all positive integers $n$ and all $m$-homogeneous polynomials
$P:\ell_{\infty}^{n}\left(  \mathbb{K}\right)  \rightarrow\mathbb{K}$, with
$\mathbb{K=R}$ or $\mathbb{C}$.

It is well-known that (\ref{comp}) is a corollary of (\ref{09}). In fact,
according to \cite[(2.1)]{fr} there are constants $C,C_{\mathbb{C}}>0$ such
that%
\begin{align}
{\sum\limits_{\left\vert \alpha\right\vert =m}}\left\vert a_{\alpha
}\right\vert  &  \leq\left(  {\sum\limits_{\left\vert \alpha\right\vert =m}%
}\left\vert 1\right\vert ^{\frac{2m}{m-1}}\right)  ^{\frac{m-1}{2m}}\left(
{\sum\limits_{\left\vert \alpha\right\vert =m}}\left\vert a_{\alpha
}\right\vert ^{\frac{2m}{m+1}}\right)  ^{\frac{m+1}{2m}}\label{0p}\\
&  \leq\left(  C^{m}\left(  1+\frac{n}{m}\right)  ^{m}\right)  ^{\frac
{m-1}{2m}}C_{\mathbb{C},BH}^{m}\left\Vert P\right\Vert _{\infty}\nonumber\\
&  \leq C_{\mathbb{C}}^{m}\cdot n^{\frac{m-1}{2}}\left\Vert P\right\Vert
_{\infty}.\nonumber
\end{align}
For related recent results we refer to \cite{Nu3, os, Se} and \cite{seip3} for
a panorama of the subject.

Both (\ref{comp}) and the polynomial (and multilinear) Bohnenblust--Hille
inequalities were originally conceived for complex scalars; the reason is that
these inequalities were motivated by problems arising over the complex scalar
field. In the last years, however, the interest in the Bohnenblust--Hille
inequality encompassed the case of real scalars, mainly due to its connections
with Quantum Information Theory (see \cite{montanaro}).

In this note we remark that (\ref{comp}) and (\ref{09}) are particular cases
of a continuum family of sharp inequalities for both complex and real scalars:

\begin{theorem}
Let $r\in\lbrack1,\frac{2m}{m+1}],$ $m,n$ be positive integers and
$\mathbb{K}=$ $\mathbb{R}$ or $\mathbb{C}.$ There is an universal constant
$L_{\mathbb{K}}>0$ such that%
\begin{equation}
\left\Vert P\right\Vert _{r}\leq L_{\mathbb{K}}^{m}\cdot n^{\left(  \frac
{m}{r}-\frac{m+1}{2}\right)  }\left\Vert P\right\Vert _{\infty}\label{ki}%
\end{equation}
for all $m$-homogeneous polynomials $P:\ell_{\infty}^{n}\left(  \mathbb{K}%
\right)  \rightarrow\mathbb{K}.$ Moreover, the power $\frac{m}{r}-\frac
{m+1}{2}$ is optimal. 
\end{theorem}

It is worth mentioning that, in general, the adaptation of asymptotic results
involving homogeneous polynomials from the complex setting to real scalars is
not a straightforward task. In fact, sometimes polynomials present a
completely different behavior when we change the scalar field from
$\mathbb{C}$ to $\mathbb{R}$ (we refer to \cite[page 58]{studia} for an
illustrative example of this fact).

\section{The proof}

Let $r\in$ $[1,\frac{2m}{m+1}].$ The proof of (\ref{ki}) for complex scalars
is easily obtained by adapting the argument used in (\ref{0p}). For real
scalars, according to \cite{CMPS}, if $P:\ell_{\infty}^{n}\left(
\mathbb{R}\right)  \rightarrow\mathbb{R}$ is an $m$-homogeneous polynomial,
then
\[
\left\Vert P_{\mathbb{C}}\right\Vert _{\infty}\leq2^{m-1}\left\Vert
P\right\Vert _{\infty},
\]
where $P_{\mathbb{C}}$ is the complexification of $P;$ this result goes back
to the Visser's paper \cite{Visser}. So, we obtain (\ref{0p}) for real
scalars. It also simple to see that the constant $L_{\mathbb{K}}$ can be
chosen independent of $m,r,n.$ Now let us prove the optimality of the exponent
$\frac{m}{r}-\frac{m+1}{2};$ for this task let us suppose that the result
holds for a power $q<\frac{m}{r}-\frac{m+1}{2}$.

For each $m,n$, let
\begin{align*}
P_{m,n}  &  :\ell_{\infty}^{n}\left(  \mathbb{K}\right)  \rightarrow
\mathbb{K}\\
P_{m,n}(w)  &  ={\sum\limits_{\left\vert \alpha\right\vert =m}}\varepsilon
_{\alpha}w^{\alpha}%
\end{align*}
be the $m$-homogeneous Bernoulli polynomial satisfying the
Kahane--Salem--Zygmund inequality (note that this inequality is also valid for
real scalars, see \cite{Nun2}).

The proof follows the lines of \cite[Theorem 10.2]{Nun2}; the essence of this
argument can be traced back to Boas' classical paper \cite{Boas}. We can
suppose $n>m.$ As in \cite{Nun2}, we have
\[
{\textstyle\sum_{\left\vert \alpha\right\vert =m}}\left\vert \varepsilon
_{\alpha}\right\vert ^{r}=p(n)+\frac{1}{m!}{\textstyle\prod\limits_{k=0}%
^{m-1}}(n-k),
\]
where $p\left(  n\right)  >0$ is a polynomial of degree $m-1.$ If (\ref{ki})
was valid with the power $q$, then there would exist a constant
$C_{q,\mathbb{K}}>0$ so that%
\begin{align*}
\left(  {\sum_{\left\vert \alpha\right\vert =m}}\left\vert \varepsilon
_{\alpha}\right\vert ^{r}\right)  ^{\frac{1}{r}}  &  \leq C_{q,\mathbb{K}}%
^{m}\cdot n^{q}\left\Vert P_{m,n}\right\Vert _{\infty}\\
&  \leq C_{q,\mathbb{K}}^{m}\cdot n^{q}\cdot C_{KSZ}\cdot n^{\left(
m+1\right)  /2}\sqrt{\log m},
\end{align*}
where $C_{KSZ}>0$ is the universal constant from the Kahane--Salem--Zygmund
inequality. Hence%
\[
C_{q,\mathbb{K}}^{m}C_{KSZ}\geq\frac{1}{n^{q}\cdot n^{\left(  m+1\right)
/2}\sqrt{\log m}}\left(  p(n)+\frac{1}{m!}{\textstyle\prod\limits_{k=0}^{m-1}%
}(n-k)\right)  ^{1/r}%
\]
for all $n$. Raising both sides to the power of $r$ and letting $n\rightarrow
\infty$ we obtain%
\[
\left(  C_{q,\mathbb{K}}^{m}C_{KSZ}\right)  ^{r}\geq\lim_{n\rightarrow\infty
}\left(  \frac{p(n)}{n^{qr}\cdot n^{r\left(  m+1\right)  /2}\left(  \sqrt{\log
m}\right)  ^{r}}+\frac{s(n)}{n^{qr}\cdot n^{r\left(  m+1\right)  /2}\left(
\sqrt{\log m}\right)  ^{r}}\right)  ,
\]
with
\[
s(n)=\frac{1}{m!}{\textstyle\prod\limits_{k=0}^{m-1}}(n-k).
\]
Since $q<\frac{m}{r}-\frac{m+1}{2}$, we have $\deg s=m>qr+r(m+1)/2$ and thus
the limit above is infinity, a contradiction.

\end{document}